\documentclass{amsart}[11pt]
\usepackage{amsmath,amssymb,pifont,calc,bm}
\usepackage{graphicx}
\usepackage{ifpdf}
\usepackage{color}
\usepackage{comment}

\newtheorem{theorem}{Theorem}[section]
\newtheorem{proposition}[theorem]{Proposition}

\newtheorem{corollary}[theorem]{Corollary}

\theoremstyle{definition}

\theoremstyle{remark}
\newtheorem{remark}[theorem]{Remark}
\newtheorem*{remarks*}{Remarks}
\newtheorem*{remark*}{Remark}

\numberwithin{equation}{section}

\newcommand{\h}{\hbar}

\newcommand{\bbR}{{\mathbb R}}

\begin{document}
\openup1pt

\title{A semiclassical heat trace expansion for the perturbed  harmonic oscillator}

\author{V. Guillemin}\address{Department of Mathematics\\
Massachusetts Institute of Technology \\Cambridge, MA 02139 \\
USA}\thanks{V. Guillemin is supported in part by NSF grant
DMS-1005696.}\email{vwg@math.mit.edu}
\author{A. Uribe}\address{Department of Mathematics \\
University of Michigan \\ Ann Arbor \\ MI \\ 48109 \\ USA}
\thanks{A. Uribe is supported in part by NSF grant DMS-0805878.}
\email{uribe@umich.edu}
\author{Z. Wang}
\address{Department of Mathematics \\
University of Michigan \\ Ann Arbor \\ MI \\ 48109 \\ USA} \email{wangzuoq@umich.edu}

\begin{abstract}
In this paper we
study the heat trace expansion of the perturbed harmonic oscillator by
adapting to the semiclassical setting techniques developed by
Hitrick-Polterovich in [HP].   We use the expansion to obtain certain
inverse spectral results.
\end{abstract}

\maketitle

\section{Introduction}

Hitrik and Polterovich obtained in \cite{HP} a simple formula for the on-diagonal heat kernel 
expansion of the Schr\"odinger operator, $-\Delta + V$, with $V \in C^\infty(\mathbb R^n)$ a bounded real-valued potential.  In this paper we apply their techniques to study the semiclassical
behavior of the on-diagonal heat kernel
expansion for the perturbed semi-classical harmonic oscillator 
\begin{equation}
\label{semiHO}
H = 
\sum_{i=1}^{n} \left( -\frac{\hbar^2}2 \frac{\partial^2}{\partial x_i^2} + \frac{x_i^2}2 - \frac {\hbar}2\right) + \hbar^2 V.
\end{equation} 
More precisely, we consider the kernel of the operator $e^{-tH}$ but with ordinary time
$t$ replaced by 
\begin{equation}
\label{heatTime}
t = \frac 1{\hbar} \log \frac{1+\hbar s}{1-\hbar s}
= 2s \left( 1+ \frac{\hbar^2 s^2}{3} + \frac{\hbar^4 s^4}{5} + \cdots \right),
\end{equation}
which greatly simplifies the calculations in the semiclassical regime.

Our first result is:

\begin{theorem}\label{Main}  Assume that $V\in C^\infty({\mathbb R}^n)$ is bounded below and 
that it and all its derivatives have at most polynomial growth at infinity.  Then, on the diagonal,
the Schwartz kernel of the operator $e^{-tH}$, where $t$ is given by (\ref{heatTime}), 
has an asymptotic expansion as $\h$ tends to zero of the form
\begin{equation}\label{mainExpansion}
\h^2\, \sum_{k=0}^\infty \h^{2k}\,\Upsilon_k(s, x).
\end{equation}
Moreover, the first three coefficients in this expansion, integrated over $\bbR^{n}$, 
determine the following quantities:  
\begin{equation}\label{mainInvariants}
\int_{\bbR^{n}} V(x)\, e^{-s|x|^{2}}\,dx,\quad\int_{\bbR^{n}} V^2(x)\, e^{-s|x|^{2}}\,dx,
\quad \int_{\bbR^{n}} \left(V^3(x)-V\Delta V\right)\, e^{-s|x|^{2}}\,dx .
\end{equation}
\end{theorem}
The Hitrik-Polterovich method, adapted to the present situation, results in a procedure to 
compute the $\Upsilon_k$.   In particular $\Upsilon_{0}(s, x) = 2sV(x) e^{-|x|^{2}s}$.

For $\h$ sufficiently small the spectrum of $H$ is discrete and the quantities (\ref{mainInvariants}) are {\em spectral invariants} of $V$ associated to the $\h$-dependent spectrum of $H$.  By analyzing these
invariants we obtain several inverse spectral results, namely:

\begin{corollary}\label{SpecResults}  Let $S_r = \{ x\in\bbR^n\;;\;|x| = r\}$.
The following properties of $V$ can be detected from the $\h$-dependent spectrum 
of $H$:
\begin{enumerate}
\item[(a)] Whether $V$ is constant on a given sphere $S_r$, and if so
the value of the constant.
\item[(b)] Whether $V$ is compactly supported, 
and if so the  the smallest annulus about the origin containing the support of $V$.
\item[(c)]  Within the class of odd functions $V$, one can determine 
whether the restriction of $V$ to any annulus about the origin is linear on that annulus.
\end{enumerate}
\begin{remark}
Item (c) is a consequence of a
much stronger but slightly more technical result (see Proposition \ref{MoreTechnical}.) 
\end{remark}

\end{corollary}

\medskip
The basic ingredients in the derivation of the expansion (\ref{mainExpansion})
are a variant of Mehler's formula and the Kantorovitz formula for expressing the heat expansion of the sum of two operators, $A$ and $B$, in terms of the heat expansion of $B$ alone.  (The latter is also the basic ingredient in the proof of the Hitrik-Polterovich result.) 
We will discuss Mehler's formula in \S 2 and the Hitrik-Polterovich formula in \S 3. 
Then in \S 4 we will describe what this formula looks like if one replaces the standard heat kernel, 
$(4\pi t)^{-\frac n2} e^{-\frac{|x-y|^2}{4t}}$, by the semi-classical Mehler kernel. 
As mentioned, the expansion of this formula in powers of $\hbar^2$ generates a sequence of heat trace invariants, 
and in \S 5 we will discuss a symbolic method for computing these invariants. 
In \S 6 we will illustrate these methods by computing the first three of these invariants, and finally in \S 7 we will 
prove the aforementioned inverse spectral results.

\section{Mehler's formula}

Let $L$ be the operator
  \begin{equation}\label{HO}
  L= \frac 12\sum  \left( -\frac{\partial^2}{\partial x_i^2} + x_i^2 -1 \right).
  \end{equation}
Mehler's formula for the Schwartz kernel of $e^{-tL}$ is
  \begin{equation}\label{Mehler}
  e^{-tL}(x,y) = \pi^{-\frac n2}(1-e^{-2t})^{-\frac n2} \exp\left\{-\frac 1{1-e^{-2t}}\left[\frac{|x|^2+|y|^2}2(1+e^{-2t})-2e^{-t}x\cdot y\right]\right\}
  \end{equation}
   (see for instance \cite{Sim} page 38).
Rescaling the variables $x$ and $y$ by the factor $1/\sqrt{\hbar}$ and $t$ by the factor $\hbar$ we get, for the heat kernel of the semi-classical harmonic oscillator 
\[
A = -\frac{\hbar^2}2 \Delta + \frac{|x|^2}2 - \frac{n\hbar}2,
\]
the expression
  \begin{equation}\label{Kernel}
  e^{-tA}(x,y) = \pi^{-\frac n2}(1-e^{-2\hbar t})^{-\frac n2} \exp\left\{-\frac 1{\hbar(1-e^{-2\hbar t})}\left[ \frac{|x|^2+|y|^2}2(1+e^{-2\hbar t})-2e^{-t\hbar}x\cdot y \right]\right\}.
  \end{equation}
The term in square brackets can be rewritten as
\[
\frac{|x|^2+|y|^2}2(1-e^{-t\hbar})^2 +e^{-t\hbar} |x-y|^2,
\]
and hence the term in curly braces is equal to  
\begin{equation}\label{expo}
-\frac{|x|^2+|y|^2}{2\hbar} \frac{1 - e^{-t\hbar}}{1+ e^{-t\hbar}} 
- \frac{e^{-t\hbar}}{\hbar (1-e^{-t\hbar})(1+e^{-t\hbar})} |x-y|^2.
\end{equation}

Now introduce the new time scale 
  \begin{equation}\label{EhTime}
  s  = \frac 1{\hbar} \frac{1-e^{-t\hbar}}{1+e^{-t\hbar}} = \frac t2\left(1 + O(t^2\hbar^2)\right)
  \end{equation}
or, alternatively,
  \begin{equation}\label{ET2}
  e^{-t\hbar} = \frac{1 - \hbar s}{1 + \hbar s}.
  \end{equation}
Then the expression (\ref{expo}) becomes
\[
- \frac{|x|^2+|y|^2}{2} s - \frac 14\left( \frac 1{\hbar^2}\frac 1{s} - s\right) |x-y|^2 
\]
or
  \begin{equation}\label{expo2}
   -\frac{|x-y|^2}{4 \hbar^2 s} - \frac{|x+y|^2}4 s ,
  \end{equation}
and hence for the Schwartz kernel of heat operator $e^{-tA}$ we get the formula
  \begin{equation}
  \label{SK}
    e^{-tA}(x,y) = 
  (4\pi \hbar)^{-\frac n2} s^{-\frac n2} (1+s\hbar)^n \exp\left( -\frac{|x-y|^2}{4 \hbar^2 s} - \frac{|x+y|^2}4 s\right).
  \end{equation}

\section{The Kantorovitz formula}

Let $A$ and $B$ be linear operators on an appropriately defined Hilbert (Banach, Frechet, $\cdots$) space which generate strongly continuous semigroups $e^{tA}$ and $e^{tB}$ 
and such that the sets of $C^\infty$ vectors satisfy: $D^\infty(A) \subset D^\infty(A+B)$. (Both conditions will be satisfied automatically in what follows.)
 Then according to Kantorovitz, \cite{Ka}, $e^{t(A+B)}$ can be expressed as a series 
  \begin{equation}\label{Ka}
 e^{t(A+B)}= (I + tX_1 + \frac{t^2}2 X_2 + \cdots) e^{tA},
  \end{equation}
where the $X_i$'s are defined by 
  \begin{equation}\label{X012}
  X_0 = I, \quad X_1 = B, \quad X_2  = B^2 + [A, B]
  \end{equation}
and in general 
  \begin{equation}
  \label{Xm}
  X_m = B X_{m-1} + [A, X_{m-1}].
  \end{equation}
There is also a simple closed form expression for $X_m$:  Letting $H = A+B$,
  \begin{equation}\label{Xm2}
  X_m = H^m  - m H^{m-1}A + {m \choose 2} H^{m-2}A^2 + \cdots.
  \end{equation}

\noindent{\underline{\bf Example:}} (\cite{HP}) Let $A = -\Delta_{\mathbb R^n}$ and $H = -\Delta_{\mathbb R^n} +V$. Then
  \begin{equation}\label{HPA}
  e^{-tA}(x,y) = (4\pi t)^{-\frac n2} e^{-\frac{|x-y|^2}{4t}}
  \end{equation}
and hence 
  \begin{equation}\label{HPH}
  e^{-tH}(x, y) = (4\pi t)^{-\frac n2} \sum_{m=0}^{\infty} (-1)^m  \frac{t^m}{m!} a_m(x, y, t),
  \end{equation}
where 
  \begin{equation}\label{HPak}
  a_m(x,y, t) = \sum_{l+j=m}{m \choose l} H_x^l \Delta_x^j e^{-\frac{|x-y|^2}{4t}}.
  \end{equation}
From this identity it is relatively easy to get an asymptotic expansion of $e^{-tH}(x,x)$ as a Taylor series in $t$ for which the summands are (at least in principle) computable. (See \cite{HP}, \S 2 for details.) 

\medskip
In the spirit of this example, let
  \begin{equation}
  \label{HO2}
  A = \sum_{i=1}^n \left( -\frac{\hbar^2}2 \frac{\partial^2}{\partial x_i^2} + \frac{x_i^2}2 - \frac {\hbar}2  \right)
  \end{equation}
and let $H = A + \hbar^2 V$. Then, as above, 
  \begin{equation}\label{HKernel}
  e^{-tH}(x,y) = \sum_{m=0}^\infty (-1)^m \frac{t^m}{m!} a_m(x, y, t, \hbar),
  \end{equation}
where 
  \begin{equation}\label{Hak}
  a_m(x, y, t, \hbar) = \sum_{l+j=m} (-1)^j {m \choose l} H_x^l A_x^j e^{-tA}(x, y).
  \end{equation}
By (\ref{SK}) the computation of this sum reduces to computing
  \begin{equation}
  \label{Target} 
  \sum_{l+j=m} (-1)^j{m \choose l} H^l A^j 
  \exp\left( -\frac{|x-y|^2}{4\hbar^2 s} \right)f(x,y, s)t^m ,
  \end{equation}
where $s$ is given by (\ref{EhTime}) and
  \begin{equation} \label{f}
  f(x,y,s) = \exp  \left(-s \frac{|x+y|^2}4\right).
  \end{equation}
The expression above is similar to the Hitrik-Polterovich expression 
  \begin{equation*}
  \sum_{l+j=m} {m \choose l}(-\Delta +V)^l \Delta^j 
  \exp\left(-\frac{|x-y|^2}{4t}\right)t^m ,
  \end{equation*}
except for the presence of the factor (\ref{f}). 
However, since we'll mainly be interested in the $\hbar$ dependence of the 
expression (\ref{Target}), and (\ref{f}) depends in an explicit way on $\hbar$, 
our computations will be very similar to theirs. 

\section{Computations}

As above let $A$ be the operator \[ \sum_{i=1}^n \left( -\frac{\hbar^2}2 \frac{\partial^2}{\partial x_i^2} + \frac{x_i^2}2 - \frac {\hbar}2  \right),\]
$B$ the operator, $\hbar^2 V$, and $X_m, m=0, 1, 2, \cdots$ the operators defined inductively by
  \begin{equation}\label{Xm3} 
  X_m = B X_{m-1} + [A, X_{m-1}]
  \end{equation}
and $X_0 = \mathrm{I}$. It will be convenient to write this formula as
  \begin{equation}\label{Xm4}
  X_m = \hbar^2 V X_{m-1} - \hbar^2[\frac {\Delta}2, X_{m-1}] + [\frac{x^2}2, X_{m-1}].
  \end{equation} 
From this formula one gets:
\begin{proposition}
The operators $X_m$ are of the form
\begin{equation}\label{XmF}
X_m = \hbar^m \sum_{\substack{i=1\\i\equiv m\text{ mod} (2)} }^{m} \hbar^i X_m^{i-1}, 
\end{equation}
where $X_m^{i-1}$ is a differential operator of degree $i-1$ not depending on $\hbar$. 
Moreover, these operators
satisfy 
\begin{equation}\label{Xmi}
X_{m+1}^i = -\left[\frac{\Delta}2, X_m^{i-1}\right]+\left[\frac{x^2}2, X_m^{i+1}\right] + VX_m^{i-1}.
\end{equation}
\end{proposition}
The proof is a simple inductive argument.

To compute the $m^{th}$ summand in the Kantorovitz expansion (\ref{Ka}), 
we must apply $X_m$ to the Mehler kernel
\begin{equation}
\label{MK}
e(x, y, s, \hbar)  = \exp{\left(-\frac{|x-y|^2}{4 \hbar^2 s} - \frac{|x+y|^2}4 s\right)}
\end{equation}
and then set $x=y$. 
\begin{proposition}
\[ 
X_m \left. \left(  e(x, y, s, \hbar)\right)\right|_{x=y}  = \hbar^m 
\sum_{\substack{i=1\\i\equiv m\text{ mod} (2)} }^m \hbar^i X_m^{i-1} 
\left( e(x, y, s, \hbar)\right)|_{x=y}
\]
is equal to:
 For $m$ odd and with $l=\frac{m-1}2$, 
\begin{equation}\label{oddm}
  \hbar^{m+1}\left(\sum_{r=0}^l e_{m,r}(x, s)\hbar^{2r}\right)s^{-l}e^{-s|x|^2}
  \end{equation}  
and for $m$ even and with $l=\frac{m}2 -1$,
  \begin{equation}\label{evenm}
  \hbar^{m+2}\left(\sum_{r=0}^l e_{m,r}(x, s)\hbar^{2r}\right)s^{-l}e^{-s|x|^2},
  \end{equation}  
where in all cases the $e_{m,r}$ are polynomials in $s$ of degree at most $2r$. 
\end{proposition}

\begin{proof}
We first note that for multi-indices, $\mu$, 
 \begin{equation}\label{der}  
  \left.\partial_x^{\mu} e^{-\frac{|x-y|^2}{4\hbar^2 s}}\right|_{x=y} = c_{\mu} \hbar^{-|\mu|} s^{-\frac{|\mu|}2}
  \end{equation}  
for even $\mu$ and $0$ for non-even $\mu$, where $c_\mu = (-\frac 14)^{|\nu|} \frac{\mu!}{\nu!}$ for $\mu=2\nu$.  The result follows from this, Leibniz' formula, and the properties of the operators $X_m^{i-1}$.

\end{proof}

Thus making the substitution 
  \begin{equation}\label{t}
 t = \frac 1{\hbar} \log \frac{1+\hbar s}{1-\hbar s} = 2s \left( 1+ \frac{\hbar^2 s^2}{3} + \frac{\hbar^4 s^4}{5} + \cdots \right),
  \end{equation}
the $t^m X_m e(x,y, t, \hbar)|_{x=y}$ term in the Kantorovitz formula gets converted into
\begin{equation}
\label{oddmexpa}
2^m \hbar^{m+1} s^{l+1} \left( \sum_{r=0}^l e_{m, r}(x, s) \hbar^{2r} \right) \left( 1+ \frac{\hbar^2 s^2}{3} + \frac{\hbar^4 s^4}{5} + \cdots  \right)^m,
\end{equation}
$l = \frac{m-1}2$, for $m$ odd, and 
\begin{equation}
\label{evenmexpa}
2^m \hbar^{m+2} s^{l+2} \left( \sum_{r=0}^l e_{m, r}(x, s) \hbar^{2r} \right) \left( 1+ \frac{\hbar^2 s^2}{3} + \frac{\hbar^4 s^4}{5} + \cdots  \right)^m,
\end{equation}
$l=\frac m2 -1$, for $m$ even.

\section{Symbolic features of the expansions (\ref{oddm})-(\ref{evenm})}

 We showed above that there exist functions $\rho_m(x,s)$, $m=0, 1,\ldots$ such that
\begin{equation}
\label{Xmeodd}
\left.X_m e(x, y, s, \hbar) \right|_{x=y} = \hbar^{m+1} \rho_m(x, s) + O(\hbar^{m+3})
\end{equation}
for $m$ odd and   
\begin{equation}
\label{Xmeeven}
\left.X_m e(x, y, s, \hbar) \right|_{x=y} = \hbar^{m+2} \rho_m(x, s) + O(\hbar^{m+4})
\end{equation}
for $m$ even. We will show in this section that for $m$ odd, $\rho_m(x, s)$ is computable purely by ``symbolic" techniques and will prove a somewhat weaker form of this assertion for $m$ even. 

Let 
\begin{equation}
\label{Xm+1i}
X_{m+1}^i = \sum_{|\alpha| \le i} a^\alpha_{i, m+1}(x) D^\alpha
\end{equation}
and let 
\begin{equation}\label{pm+1i}
p_{m+1}^i  = \sum_{|\alpha| \le i} a^\alpha_{i, m+1}(x) \xi^\alpha
\end{equation} 
be the full symbol of $X_{m+1}^i$
From (\ref{Xmi}) and standard composition formula for left Kohn-Nirenberg symbols one gets 
\begin{equation}\label{symbol}
p^i_{m+1}   =  \left(\sum_{r=1}^{n} \frac {\xi_r}{\sqrt{-1}}  \frac{\partial}{\partial x_r} - \frac 12 \frac{\partial^2}{\partial x_r ^2}\right) p_m^{i-1} 
 + \left(\sum_{r=1}^{n} {\sqrt{-1}} x_r \frac{\partial}{\partial \xi_r} + \frac 12 \frac{\partial^2}{\partial \xi_r ^2}\right) p_m^{i+1} 
+ V(x)p_m^{i-1}.
\end{equation}
In particular, if 
\[\sigma_{m+1}^i = \sum_{|\alpha|=i} a^\alpha_{i, m+1}(x) \xi^\alpha\]
is the principal symbol of $X^i_{m+1}$ and 
\[
\tilde \sigma^{i}_{m+1} =   \sum_{|\alpha|=i-1} a^\alpha_{i, m+1}(x) \xi^\alpha
\]
the subprincipal symbol, we get from (\ref{symbol}) that
\begin{equation}
\label{sigmaim+1}
\sigma^i_{m+1} = \frac 1{\sqrt{-1}} \sum_{r=1}^{n} \left(\xi_r \frac{\partial}{\partial x_r} \sigma_m^{i-1} - x_r \frac{\partial}{\partial \xi_r} \sigma_m^{i+1}\right)
\end{equation}
and 
\begin{equation}
\label{tsigmaim+1}
\tilde \sigma^i_{m+1} = \frac 1{\sqrt{-1}} \sum_{r=1}^{n} \left( \xi_r \frac{\partial}{\partial x_r} \tilde \sigma_m^{i-1} - x_r \frac{\partial}{\partial \xi_r} \tilde \sigma_m^{i+1}\right) - \frac 12 \sum_{r=1}^{n} \left( \frac{\partial^2}{\partial x_r ^2} \sigma_m^{i-1} -\frac{\partial^2}{\partial \xi_r ^2} \sigma_m^{i+1} \right)+V \sigma_m^{i-1}.
\end{equation}
Letting 
\[\sigma_m = \hbar^m \sum_{\substack{i=1\\i\equiv m\text{ mod} (2)} }^{m} \hbar^i \sigma_m^{i-1}, 
\] 
and letting $\mathcal U$ be the raising operator 
\[
\mathcal U \hbar^i \sigma = \hbar^{i+2}\sigma
\]
for $i \ge 0$, we can write these formulas more succinctly in the form
\begin{equation}\label{sigmam}
\sigma_m = \frac 1{\sqrt{-1}}  \sum_{r=1}^n \left(\xi_r \frac{\partial}{\partial x_r} \mathcal U - x_r \frac{\partial}{\partial \xi_r}  \right) \sigma_{m-1}
\end{equation}
and 
\begin{equation}\label{tsigmam}
\tilde \sigma_m = \frac 1{\sqrt{-1}}  \sum_{r=1}^n \left(\xi_r \frac{\partial}{\partial x_r} \mathcal U - x_r \frac{\partial}{\partial \xi_r}  \right) \tilde \sigma_{m-1} + \sum_{r=1}^n \left( -\frac 12 \frac{\partial^2}{\partial x_r^2} \mathcal U + \frac 12 \frac{\partial^2}{\partial \xi_r^2} \right) \sigma_{m-1} + V \mathcal U \sigma_{m-1}.
\end{equation}
In particular, iterating (\ref{sigmam}) we get 
\begin{equation}\label{sigmam2}
\sigma_m = \left[ \frac 1{\sqrt{-1}} \sum_{r=1}^{n} \left( \xi_r \frac{\partial}{\partial x_r} \mathcal U - x_r \frac{\partial}{\partial \xi_r} \right) \right]^{m-1} \hbar^2 V
\end{equation}
and, as special cases of (\ref{sigmam2}),
\begin{equation}\label{sigmamm-1}
\sigma_m^{m-1} = \left( \frac 1{\sqrt{-1}} \sum_{r=1}^{n} \xi_r \frac{\partial}{\partial x_r}\right)^{m-1}V.
\end{equation}

As applications of these formulas let $m$ be odd and consider the $i$th summand of
\[
\left.X_m e(x, y, s, \hbar)\right|_{x=y} =\left.\hbar^m 
\sum_{\substack{i=1\\i\equiv m\text{ mod} (2)} }^{m}
 \hbar^i X_m^{i-1} e(x, y, s, \hbar) \right|_{x=y}.
\]
By (\ref{Xm+1i}) this is equal to 
\begin{equation}\label{mid}
\left.\hbar^{m+i} \sum_{|\alpha| = i-1} a^\alpha_{i-1, m}(x)
D^\alpha_{x} e^{-\frac{|x-y|^2}{4\hbar^2 s}} \right|_{x=y} e^{-s|x|^2}
\end{equation}
plus terms of order $O(\hbar^{m+3})$ and by (\ref{der}), (\ref{mid}) is equal to 
\begin{equation}
\hbar^{m+1} \left( \sum_{|\alpha|=i-1} a^\alpha_{i-1, m}(x) c_\alpha \right) s^{-\frac{i-1}{2}} e^{-s|x|^2}.
\end{equation}
To summarize:
\begin{proposition}  For $m$ odd, the leading term $\rho_m$ in (\ref{Xmeodd}) is given by
\[
\rho_m(x, s) = e^{-s|x|^2} \sum_{\substack{i=1\\i\equiv m\text{ mod} (2)} }^{m}
s^{\frac{1-i}{2}}\sum_{|\alpha|=i-1} a^\alpha_{i-1, m}(x) c_\alpha.
\]
For each $i$, the quantity $\sum_{|\alpha|=i-1} a^\alpha_{i-1, m}(x) c_\alpha$ is obtained from the
principal symbol of $X_m^{i-1}$ by substituting every monomial $\xi^\alpha$ by the constant
$c_\alpha$.
\end{proposition}

For $m$ even the computation above is similar, however one gets $\hbar^{m+2}$ contributions to (\ref{Xmeeven}) from both the terms 
\[
\left.\hbar^{m+i} \sum_{|\alpha|=i-1} a^\alpha_{i-1, m}(x) D^\alpha e^{-\frac{|x-y|^2}{4\hbar^2 s} } e^{-\frac s4 |x+y|^2} \right|_{x=y}
\]
\underline{and} the terms 
\[  
\left.\hbar^{m+i} \sum_{|\alpha|=i-2}  a^\alpha_{i-1,m}(x) D^\alpha e^{-\frac{|x-y|^2}{4\hbar^2 s}}
 \right|_{x=y} e^{-s|x|^2}.
\]
The second summand (involving the subprincipal symbol of $X_m^{i-1}$) is as before,
\begin{equation}\label{2ndsum}
\hbar^{m+2} \sum_{|\alpha|=i-2} a^\alpha_{i-1, m}(x) c_\alpha s^{-\frac{i-2}2} e^{-s|x|^2},
\end{equation}
but the first summand (involving the principal symbol of $X_m^{i-1}$) becomes 
\begin{equation}\label{1stsum}
\hbar^{m+2} \sum_{\substack{|\alpha|=i-1 \\1\leq r\leq n}} a^\alpha_{i-1, m}(x) c_{\alpha^{(r)}} s^{-\frac{i-4}2} x_r e^{-s|x|^2},
 \end{equation}
where $\alpha^{(r)} = (\alpha_1, \cdots, \alpha_r-1, \cdots, \alpha_n)$. This proves:

\begin{proposition}
For $m$ even the leading term of (\ref{Xmeeven}) depends only on the principal and subprincipal symbols of $X_m$. 
\end{proposition}
 
 \medskip
 We now explore some spectral consequences of the previous results.
 \begin{proposition}
 For $m$ odd the quantities $\int \rho_m(x,s)\, dx$ are spectral invariants of $V$.
 \end{proposition}
 \begin{proof}
For $m$ odd we can, by (\ref{sigmam2}), express $\rho_m(x, s)$ as a sum of terms of the form 
$x^\alpha \frac{\partial^\beta V}{\partial x^\beta} e^{-s|x|^2}$, where $|\alpha|+|\beta| \le m-1$.  The associated contribution to the heat trace 
\begin{equation}\label{heatint}
\int x^\alpha \frac{\partial^\beta V}{\partial x^\beta} e^{-s|x|^2} dx
\end{equation}
can, by integration by parts, be written as sums of integrals of the form 
\[
\int x^\gamma V e^{-s|x|^2} dx, \qquad |\gamma|\le m-1.
\]
Thus 
\begin{equation}\label{invariant0}
\int \rho_m(x, s) e^{-s|x|^2} dx  = \int p(x, s) V e^{-s|x|^2} dx,
\end{equation}
where $p(x, s)$ is a {\em universal} polynomial of degree $m-1$ in $x$. Moreover, for every $A \in SO(n)$ the heat trace expansion for the potentials $V$ and $V^A$, where $V^A(x) = V(Ax)$, are the same. Hence by averaging over $SO(n)$ we can assume that $p(x, s)$ is $SO(n)$ invariant, i.e. 
\[
p(x, s) = \sum_{i=0}^k \chi_i(s) |x|^{2i}, \qquad k=\frac{m-1}2
\]
and thus (\ref{invariant0}) becomes
\begin{equation}\label{invariant}
\int \rho_m(x, s) e^{-s|x|^2} dx = \int_0^\infty dr \sum_{i=1}^k \chi_i(s) r^{2i} e^{-sr^2} \int_{|x|=r} V(x) d\sigma_r,
\end{equation}
where $d\sigma_r$ is the standard volume form on the $(n-1)$-sphere $|x|=r$. We will see below however that, for each $r>0$,  the integral 
\begin{equation}
\int_{|x|=r} V(x) d\sigma_r
\end{equation}
is itself a spectral invariant of the perturbed harmonic oscillator and hence the terms $\hbar^{m+1} \int \rho_m(x, s) dx$ in the heat trace expansion above can be read off from it. 
\end{proof}

In the case $m$ even one also gets a similar description of the contributions of (\ref{Xmeeven}) to the heat trace. The contribution coming  from the term (\ref{1stsum}) only depends on $\sigma_m^{i-1}$ and hence as above is expressible in terms of (\ref{invariant}); and as for the contributions coming from (\ref{2ndsum}) one can prove by induction that these give rise to heat trace invariants which are sums of expressions of the form (\ref{heatint}) and of the form
\begin{equation}\label{heatint2}
\int x^\alpha V \frac{\partial^\beta V}{\partial x^\beta} e^{-s|x|^2}\ dx.
\end{equation}
Indeed the second summand in (\ref{tsigmaim+1}) is purely symbolic; so as we've just seen it contributes terms of type (\ref{heatint}) to the heat trace. Similarly the third summand contributes terms of type (\ref{heatint2}) and by a simple induction on $m$ one can show that the first summand of (\ref{tsigmaim+1}) is a linear combination of terms of the form 
\[
\xi^\alpha x^\beta \frac{\partial^\gamma V}{\partial x^\gamma }\frac{\partial^\delta V}{\partial x^\delta}
\]
and 
\[
\xi^\alpha x^\beta \frac{\partial^\gamma V}{\partial x^\gamma}.
\]
These give rise to  contributions to the heat trace of the form (\ref{heatint}) and 
\[
\int x^\mu \frac{\partial^\nu V}{\partial x^\nu} \frac{\partial^\gamma V}{\partial x^\gamma} e^{-s|x|^2} dx,
\]
which by integration by parts can be written as expressions of the form (\ref{heatint2}). Finally, the $O(n)$ invariance of the heat trace enables one to simplify these further and rewrite them as sums of the form 
\[
\int |x|^{2j} V(\sum x_i \frac{\partial}{\partial x_i})^k \Delta^l V dx.
\]

\section{The first heat invariants}

It is easy to see (either by direct computation or by the symbolic formulas in the preceding section) that for $m \le 4$ the $X_m^{i-1}$'s are given by 
\begin{equation}
\label{X10}
X_1^0 = V,
\end{equation}
\begin{equation}
\label{X21}
X_2^1 = -\sum \frac{\partial V}{\partial x_i} \frac{\partial}{\partial x_i} + V^2 - \frac{\Delta V}2, 
\end{equation}
\begin{equation}
\label{X32}
X_3^2 = \sum \frac{\partial^2 V}{\partial x_i \partial x_j} \frac{\partial^2\  }{\partial x_i \partial x_j}  + \sum \frac{\partial}{\partial x_i}(\Delta V - \frac 32 V^2)\frac{\partial}{\partial x_i} + \frac{\Delta^2}4 V -\frac{\Delta}2 V^2 + V^3 -\frac{V\Delta V}2,
\end{equation}
\begin{equation}
\label{X30}
X_3^0 = \sum x_i \frac{\partial V}{\partial x_i},
\end{equation}
\begin{equation}
\label{X43}
X_4^3 = [-\frac {\Delta}2, \sum  \frac{\partial^2 V}{\partial x_i \partial x_j} \frac{\partial^2\  }{\partial x_i \partial x_j}  + \sum \frac{\partial}{\partial x_i}(\Delta V - \frac 32 V^2)\frac{\partial}{\partial x_i}] + 
\sum V \frac{\partial^2 V}{\partial x_i \partial x_j} \frac{\partial^2\  }{\partial x_i \partial x_j} 
\end{equation}
plus terms of degree less than two, and 
\begin{equation}\label{X41}
X_4^1 = -\left[\frac{\Delta}2, \sum x_i \frac{\partial V}{\partial x_i}\right] + 
\left[\frac {x^2}2,  \sum \frac{\partial^2 V}{\partial x_i \partial x_j} \frac{\partial^2\  }{\partial x_i \partial x_j}  + \sum \frac{\partial}{\partial x_i}(\Delta V - \frac 32 V^2)\frac{\partial}{\partial x_i}\right] +\frac 12 \sum x_i \frac{\partial}{\partial x_i} V^2.
\end{equation}

Thus the $\hbar^2$ term in the heat trace expansion determines 
\begin{equation}
\int V e^{-s|x|^2} dx,
\end{equation} 
and hence by the inverse Laplace transform determines the integral 
\begin{equation}\label{1stinv}
\int_{|x|=r} Vd\sigma_r.
\end{equation}
for each $r >0$.

The $\hbar^4$ term involves the $m=3$ contribution of $(\ref{Xmeodd})$, but as we saw above this is expressible in terms of (\ref{1stinv}). As for the contribution of (\ref{Xmeeven}) to the $\hbar^4$ term, the first and third summands can be converted by integration by parts into integrals which are expressible in terms of (\ref{1stinv}) and the second summand gives a new heat invariant, 
\begin{equation}\label{2ndinvg}
\int V^2 e^{-s|x|^2} dx,
\end{equation}
which, by the inverse Laplace transform, is convertible into 
\begin{equation}\label{2ndinv}
\int_{|x|=r} V^2 d\sigma_r.
\end{equation}

The $\hbar^6$ term in the heat trace expansion involves the $m=5$ contribution of (\ref{Xmeodd}) which, as we saw in the previous section, is expressible in terms of (\ref{1stinv}), the $m=4$ contribution of (\ref{Xmeeven}), which is subprincipal and hence only involves the terms in (\ref{X41}) (all of which can be converted, by integration by parts, into expressions in (\ref{1stinv}) and (\ref{2ndinv})), and the cubic and quadratic terms in (\ref{X43}) all of which, except for the term
\begin{equation}\label{VV}
\sum V \frac{\partial^2 V}{\partial x_i \partial x_j}\frac{\partial}{\partial x_i}\frac{\partial}{\partial x_j},
\end{equation}
can be converted by integration by parts into expressions in (\ref{1stinv}) and (\ref{2ndinv}). Finally the $\hbar^6$ terms coming from (\ref{X32}) and (\ref{X30}) are all convertible by integration by parts into expressions in (\ref{1stinv}) and (\ref{2ndinv}) except for the last summand of (\ref{X32}): the term 
\begin{equation}\label{V3}
V^3 - \frac{V\Delta V}2.
\end{equation}
The term (\ref{VV}) gives, by (\ref{2ndsum}) and (\ref{evenmexpa}) a contribution 
\begin{equation}
-\hbar^6 s^2 \frac{V\Delta V}{2} e^{-s|x|^2}
\end{equation}
to the heat trace expansion, and the term (\ref{V3})  gives, by (\ref{oddmexpa}), a contribution 
\begin{equation}
\hbar^6 s^2 (V^3 - \frac{V\Delta V}2 ) e^{-s|x|^2}
\end{equation}
to the heat trace expansion; hence the sum of these two terms gives rise to a new heat trace invariant
\begin{equation}\label{3rdinvg}
\int (V^3 - V\Delta V) e^{-s|x|^2} dx
\end{equation}
which, by the inverse Laplace transform, can be converted into the invariant
\begin{equation}\label{3rdinv}
\int_{|x|=r} (V^3 - V\Delta V)d\sigma_r.
\end{equation}
This finishes the proof of Theorem \ref{Main}.

\section{Applications to inverse spectral problems}

In this section we apply the first heat invariants (\ref{1stinv}), (\ref{2ndinv}) and  
(\ref{3rdinv}) above to the inverse spectral problems of recovering information about
$V$ from the $\h$-dependent spectrum of $H$.

Fix any $r>0$. First let's consider $V$ which minimize the second invariant (\ref{2ndinv}) subject to the constraint
\begin{equation}
\label{constraint}
\int_{|x|=r} V d\sigma_r = \mathrm{constant}.
\end{equation}
According to the Cauchy-Schwartz inequality, the minimizers are exactly those functions $V$ that are constant on the sphere $|x|=r$. It follows that the set of potentials $V$ that are constant on a given sphere $|x|=r$ is intrinsically defined by its spectral properties. Moreover, one can spectrally determine the constant value for each potential in this set.   This proves parts
(a) and (b) of Corollary \ref{SpecResults}.

\newcommand{\Sr}{L^2(S_r)}
Next let's assume that $V$ is an odd potential. (Using band invariant techniques, one can show that being odd is also a spectral property, c.f. \cite{GUW}.)   As we have seen that the invariants determines $\|V\|_{\Sr}^2$ for each $r$.  Taking the $r$ derivative of the second invariant, we get that
 $\langle V, \frac{\partial  V}{\partial r}\rangle_{\Sr}$ is a spectral invariant. On the other hand, since $V$ is odd, 
the third invariant (\ref{3rdinvg}) becomes 
\begin{equation}\label{3rdinvodd}
- \int V \Delta V e^{-sr^2} r^{n-1}drd\sigma_r. 
\end{equation}
Recall that in spherical coordinates
\[
\Delta V = \frac{\partial^2 V}{\partial r^2} + \frac{ n-1}{r} \frac{\partial V}{\partial r} + \frac 1{r^2} \Delta_{S_r}V.
\]
A simple computation shows that
\[
- \int V \frac{n-1}r \frac{\partial V}{\partial r} e^{-sr^2} r^{n-1}drd\sigma_r = \frac{n-1}2 \int V^2 \frac{d}{dr}\left(  {e^{-sr^2}}{r^{n-2}} \right) drd\sigma_r
\]
which is a spectrally determined quantity, since we know the integrals (\ref{2ndinv}) for all $r$. Similarly,
\[
-\int V \frac{\partial^2 V}{\partial r^2} e^{-sr^2} r^{n-1} drd\sigma_r = \int \left(\frac{\partial V}{\partial r}\right)^2 e^{-sr^2}r^{n-1} drd\sigma_r  + \int V \frac{\partial V}{\partial r} \frac{d}{dr}(e^{-sr^2}r^{n-1}) drd\sigma_r ,
\]
and again the second term on the right
\[
 \int V \frac{\partial V}{\partial r} \frac{d}{dr}(e^{-sr^2}r^{n-1}) drd\sigma_r 
\]
is also spectrally determined according to (\ref{2ndinv}). It follows from (\ref{3rdinvodd}) 
that the integral
\begin{equation}\label{3rdinvodd2}
\int  \left(\frac{\partial V}{\partial r}\right)^2 e^{-sr^2} r^{n-1} drd\sigma_r  - \int \frac  1{r^2}V \Delta_{S_r}V e^{-sr^2} r^{n-1} drd\sigma_r
\end{equation}
is spectrally determined, which, by the inverse Laplace transform, gets converted to the invariant
$\|\frac{\partial V}{\partial r}\|_{\Sr}^2 +\frac 1{r} \langle V, -\Delta_{S_r} V\rangle_{\Sr}$ for each $r>0$. Note 
 next that for every $V$ one has the following inequality
\[
\|V\|_{\Sr}^2 \left(\|\frac{\partial V}{\partial r}\|_{\Sr}^2 + \frac{1}{r}\,\langle V, -\Delta_{\Sr} V\rangle_{\Sr} \right) \ge \langle V, \frac{\partial V}{\partial r} \rangle_{\Sr}^2 + \frac{\lambda_1}{r}\, \|V\|_{\Sr}^4,
\]
where $\lambda_1$ is the first eigenvalue of the (non-negative) Laplacian on $S_r$.  Both sides 
of the inequality are spectral invariants, and equality holds if and only if $V$ satisfies the conditions 
\begin{equation}\label{zwSpectralProperty}
\frac{\partial V}{\partial r}|_{S_r} = \chi V|_{S_r}\quad\text{and}\quad V|_{S_r}\ 
\text{is a spherical harmonic of degree one}
\end{equation}
where $\chi$ is a constant, and if so one can determine $\chi$.  This proves
the first part of the following
\begin{proposition}\label{MoreTechnical}
The class of functions defined by the conditions (\ref{zwSpectralProperty}) is spectrally
determined.
Moreover, for any $V$ in this class one can determine the ratio 
\[
\chi = \frac{\partial V}{\partial r} / V
\]
on a given sphere $S_r$.
\end{proposition}
The determination of $\chi$ (which of course can depend on $r$) is done 
by looking at the quotient of $\langle V, \frac{\partial V}{\partial r}\rangle_{L^{2}(S_{r})} / \|V\|_{L^{2}(S_{r})}^2$. 

As a consequence, one can determine whether a potential is of the form
\[
V(x) =  f(r)g(\sigma)
\]
on a given annulus $r_1 \le |r| \le r_2$, where $g(\sigma)$ is a spherical harmonic of degree one, and if so, determine the function $f(r)$.
In particular, one can spectrally determine linear potentials on any annular region $r_1 \le |r| \le r_2$: They are just the potentials in the previous class with $\chi = r$.


\end{document}